\documentclass[11pt]{article}
\usepackage[letterpaper, left=1.2in, right=1.2in, bottom=1.2in, top=1.2in]{geometry}

\usepackage[utf8]{inputenc}

\usepackage{amsmath}
\usepackage{amsfonts}
\usepackage{graphicx}
\usepackage{epstopdf}
\usepackage{color}
\usepackage{float}
\usepackage{comment}
\usepackage{standalone}
\usepackage{tabularx}

\usepackage{subcaption}
\usepackage{hyperref}

\DeclareMathOperator{\diag}{diag}

\usepackage{pgf}
\usepackage{pgfplots,tikz}
\pgfplotsset{compat=newest}
\usetikzlibrary{plotmarks}
\usetikzlibrary{arrows.meta}
\usepgfplotslibrary{patchplots}
\usepackage{grffile}
\usepackage{microtype}

\usetikzlibrary{intersections,calc,patterns,angles,quotes}
\usetikzlibrary{decorations.shapes,arrows,calc,decorations.markings,shapes,decorations.pathreplacing,shapes.arrows,arrows.meta	}
\usepackage{tikz-3dplot}
\usetikzlibrary{external}
\newcommand{\R}{\mathbb R}
\newcommand{\xin}{x_{\mathrm{init}}}
\newcommand{\xfin}{x_{\mathrm{fin}}}

\newcommand{\U}{\mathcal U}

\newcommand{\cA}{\mathcal A}
\newcommand{\X}{\mathcal X}
\newcommand{\z}{\mathcal Z}

\usepackage[nodisplayskipstretch]{setspace}

\title{\LARGE \bf
Geometric Heat Flow Method for Legged Locomotion Planning
}
\author{Yinai Fan, \and Shenyu Liu, \and Mohamed-Ali Belabbas}

\begin{document}

\maketitle

\begin{abstract}
We propose in this paper a motion planning method for legged robot locomotion based on Geometric Heat Flow framework. The motion planning task is challenging due to the hybrid nature of dynamics and contact constraints. We encode the hybrid dynamics and constraints into Riemannian inner product, and this inner product is defined so that short curves correspond to admissible motions for the system. We rely on the affine geometric heat flow to deform an arbitrary path connecting the desired initial and final states to this admissible motion. The method is able to automatically find the trajectory of robot's center of mass, feet contact positions and forces on uneven terrain.
\end{abstract}

\section{Introduction}

Planning dynamic motions of legged robots has become an increasingly important topic, due in part to improved  robot design and hardware, and in part to higher on-board computational capacity.  Typical examples of such robots designs include the MIT Cheetah~\cite{6880316}, the bipedal robot Cassie made by Agility Robotics and the Salto robot~\cite{Haldaneeaag2048}. The major  difficulty for legged locomotion planning lies in its {\it hybrid nature}: the dynamics of legged robots is governed by a  set of equations and constraints depending on whether there is contact with the ground. Hybrid systems are well-known to be difficult to handle; in fact, open questions remain even in the case of linear dynamics~\cite{liberzon2003switching}. {This difficulty stems in part from the fact that most motion planning methods do not admit natural or obvious extensions to the handle hybrid dynamics, often resulting in ad-hoc modifications that are difficult to analyze.} 
In this paper, we show that the geometric method we proposed for motion planning  extends naturally to handle hybrid dynamics. Precisely, we show how the Ansatz developed in~\cite{7963599,AGHF2019, IFAC_Mechatronics2019}, contending  that motion planning problems can be encoded into Riemannian metrics, can be applied  and using fast solvers for parabolic partial differential equations (PDE), we can obtain efficiently natural motions. 

A  variety of dynamic models have been used for different types of locomotion. Among which, the Linear Inverted Pendulum model is a simplified model that is used in cooperation with Zero Moment Point as stability criterion for bipedal walking \cite{1241826}. At the opposite extreme, full dynamics are used in order to plan the joint trajectories for motions with external contacts in~\cite{posa2014direct,neunert2017trajectory}. Centroidal dynamics, i.e. the dynamics of robot projected at its Center of Mass (CoM) \cite{orin2013centroidal}, has been used for generating whole body motions of humanoid robot \cite{dai2014whole} and hydraulic quadruped robot locomotion \cite{aceituno2017simultaneous}. By using legs with light weight, or assuming the legs do not significantly deviate from their nominal pose, one can simplify the centroidal dynamics to single rigid body dynamics 
For example, high speed bounding motions are achieved for quadruped~\cite{7139918}, biped and quadruped locomotion in complex terrains are planned in~\cite{winkler2018gait}. In this work, we use the 2 dimensional model with massless legs, which we called a {\it \textbf{Single Rigid Body Model}}.

Trajectory optimization (TO) formulates  motion planning  into an optimization problem and is widely used in legged locomotion planning, e.g. in direct collocation and differential dynamic programming methods. In \cite{dai2014whole,posa2014direct}, the hybrid dynamics of contact is modeled as complementarity problem, with the ability to plan the contact locations. By convex modeling of the dynamics, convex optimization techniques are used in~\cite{dai2016planning,ponton2016convex} for faster convergence, while the footholds need be pre-planned. The discrete nature of contact can be also modeled utilizing binary valued decision variables and solved by mix-integer solver, \cite{aceituno2017simultaneous}. More recent work~\cite{winkler2018gait} plans both gait timings and contact locations automatically by modeling the contact dynamics individually.
On the other hand, motion planning using {\it \textbf{Geometric Heat Flow (GHF)}} originates from the field of Riemannian Geometric Analysis \cite{Jost2011}. 
The method consists of encoding  the dynamic constraints into a Riemannian inner product and reduces  motion planning  to a curve shortening problem. The algorithm starts from an arbitrary path between the desired initial and final states and ``deforms'' it into an almost feasible trajectory of the system. In~\cite{7963599, LBMotionSktechArxiv2019}, fundamental convergence theorems and algorithms are proposed for driftless control affine system.  For a driftless system, The convergence to a feasible motion is guaranteed for any initial guess, under some proper assumptions. 
The evolution of trajectory from arbitrary curve to planned motion is given by the {\it geometric heat flow (GHF)}, whose solution is obtained via a PDE solver, instead of optimization solvers used in TO. 
These methods  were extended in~\cite{AGHF2019}, in which the {\it \textbf{Affine Geometric Heat Flow (AGHF)}} is introduced to handle systems with drift. Motions for robot systems are generated using AGHF in \cite{WROCO2019,IFAC_Mechatronics2019}. The present work is the first attempt at formulating the hybrid legged locomotion problem into AGHF framework and shows the potential to extend to more complicated systems and scenarios. The proposed algorithm can  encode variety of constraints including contact constraints.
Meanwhile, the convergence is guaranteed at given order, as explained in \cite{AGHF2019}.

 The remainder of the paper is organized as follows. In Section \ref{HF}, preliminary background for motion planning with GHF is introduced. Section \ref{dynamics} presents the robot dynamics and constraints for legged robot locomotion. The core of this work lies in Section \ref{forlumation}, which is the formulation of the legged locomotion into GHF frame. Section \ref{applications} shows the performance of the proposed method by examples with different scenarios. 
\section{Motion Planning using the Affine Geometric Heat Flow}\label{HF}

\subsection{Motion Planning Problem}

Consider a {\it controllable} system which is {\it affine in the control}:
\begin{equation}\label{eq:sys}
    \dot x=F_d(x)+F(x)u
\end{equation}
where $x \in \R^n$, $ u \in \R^m$, $F_d(x)$ is the drift or uncontrolled  dynamics and the columns of $F(x)\in\R^{n\times m}$ are the actuated directions of motion. Both $F_d(x), F(x)$ are assumed to be at least $C^2$ and Lipschitz. The system is {\it under-actuated} if $m<n$. We assume that $F(x)$ is of {\it constant} column rank almost everywhere in $\R^n$.

Denote by $\xin$ and $\xfin$  the desired initial and final states respectively, and by $T >0$ the time allowed to perform the motion. The space of {\it continuous controls} is denoted by $\U:=C^0([0,T]\to\R^m)$, and the space of  differentiable curves joining desired initial and final states is denoted by $\X:=\{x(\cdot)\in C^1([0,T]\to\R^n):x(0)=\xin,x(T)=\xfin\}$. We call any $x(\cdot)\in\X$ an {\it admissible curve} if there exists $u \in \U$ so that the generalized derivative $\dot x(t)$ of $x(\cdot)$ at $t$ satisfies~\eqref{eq:sys}~\cite{royden2010real}. Denote by $\X^*\subseteq\X$ the set of admissible curves. The motion planning problem is {\it feasible} if $\X^*\neq\emptyset$.

\subsection{A Brief Overview of the Affine Geometric Heat Flow}\label{AGHF}

\begin{figure}
\centering
\begin{subfigure}[t]{0.55\columnwidth}
\centering
\includegraphics{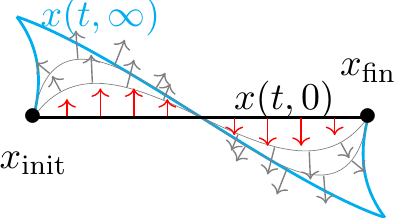}
\caption{\footnotesize Homotopy of trajectories for unicycle parallel parking.}
\label{fig:AGHF}
\end{subfigure}
\begin{subfigure}[t]{0.4\columnwidth}
\centering
\scalebox{0.45}{
\includegraphics{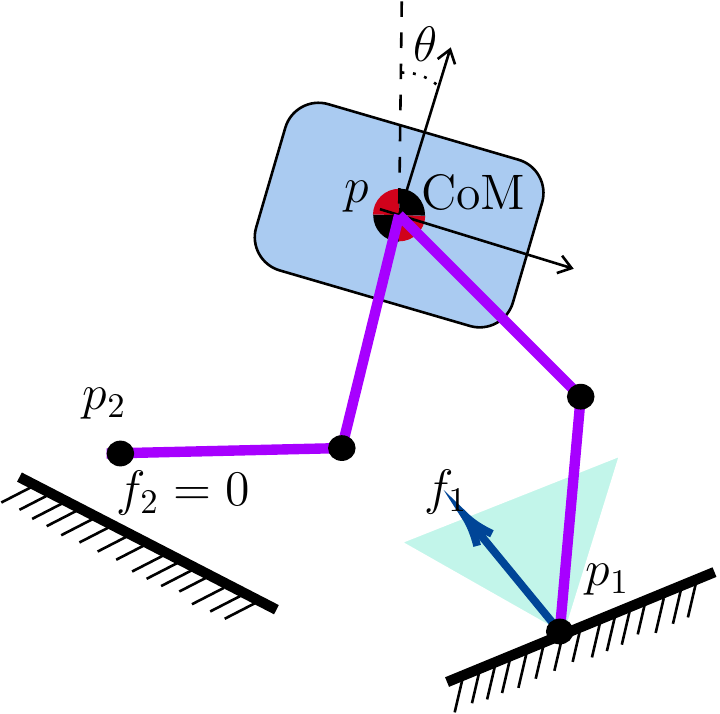}
}
\caption{\footnotesize 2D Single Rigid Body Model}
\label{fig:massless_leg_model}
\end{subfigure}

\end{figure}

The core of the planning algorithm is the {\it affine geometric heat flow (AGHF)}~\cite{AGHF2019}. This flow starts from a curve $x(t,0)$ joining $\xin$ to $\xfin$, shown by black line in Fig.~\ref{fig:AGHF}, this curve can essentially be chosen arbitrary and, in particular, does not meet the constraints, dynamic, non-holonomic or holonomic. It then ``deforms'' this curve into an admissible curve or, precisely, finds a one-parameter family of curves $x(t,s): [0,T]\times[0,s_{\max})$ where for {\it each $s$ fixed } $x(t,s), t \in[0,T]$ is a curve joining $\xin=x(0,s)$ to $\xfin=x(T,s)$, and $x(\cdot,\infty)$ is a trajectory meeting the design requirements, shown by cyan line.

To obtain this transformation from arbitrary curve to feasible trajectory, the key step is to define a Riemannian metric which  encodes the various constraints on the motion planning problem.%{\color{blue}The main insight is the following:  }
A  Riemannian metric allows us to define the {\it length} of trajectories in $\X$ and paths of shortest length %{\color{blue}, or geodesics,} 
have been widely studied; %{\color{blue}; leveraging this fact}, 
our approach is to construct a metric so that trajectories of {\it smallest length} between %{\color{blue} two points} 
$\xin$ and $\xfin$ are {\it admissible} trajectories for the system, provided such trajectories exist. We then proceed to find such curves of minimal length via a homotopy initialized at an {\it arbitrary} curve joining $\xin$ to $\xfin$. We refer the reader to~\cite{AGHF2019} for a detailed presentation.

 The interpretation of the method just given addresses systems {\it without} drift, for which a clear parallel can be made between trajectory length and admissible paths. 
 As just mentioned, a key step of the method is to define an appropriate positive definite matrix $G(x)$, which we refer to as the Riemannian {\it metric tensor}, or simply Riemnannian metric. 

In order to handle dynamics with drift, we introduced  the {\it actuated length} of a curve.  More precisely, whereas the length is given by the integral $\int_0^T ({\dot x^\top G(x) \dot x})^{1/2} dt$,  the \emph{actuated length} of a curve $x(\cdot)$  is given by 
\begin{equation}\label{actuated_length}
    \cA(x(\cdot)):=\int_0^T \underbrace{\left({(\dot x-F_d)^\top G(x) (\dot x-F_d)}\right)^{1/2}}_{L(x,\dot x)} dt.
\end{equation}
where $G(x)$ is the inner product matrix in the state space,   we will describe it below.

In the case of systems without drift, we rely on the geometric heat flow equation~\cite{Jost2011} to obtain a curve of minimal length. We introduced the AGHF to deal with actuated length. Similar to its counterpart the GHF, it is a parabolic PDE which evolves a curve with fixed end-points toward a curve of {\it minimal  actuated length}, i.e., a minimizer of~\eqref{actuated_length}. We refer to~\cite{AGHF2019} for a precise description and  intuitive explanation of the AGHF, as well as convergence results . The flow has the general form
{
\begin{equation}\label{eq:AGHF_general}
\frac{\partial x(t,s)}{\partial s} = \Psi(x(t,s),\dot x(t,s), t)
\end{equation}
where the AGHF $\Psi(x,\dot x,t)$ does not contain any differentials of $x$ with respect to $s$. Hence, one could think of the PDE as describing the deformation (omitting the boundary conditions) of a curve according to $x(t,s+\delta s) = x(t,s) + \delta s \Psi(x(t,s),\dot x(t,s), t)$, for some small step-size $\delta s$, for all $t \in [0,T]$. Fig.~\ref{fig:AGHF} shows the deformation of a parallel parking motion of a unicycle. The initial curve $x(t,0)$ is a straight line connecting initial and final position, which is not  feasible  since it requires the unicycle to side slip. The curve is deformed by the AGHF (red arrows) and eventually reach the ``Z'' shaped  curve $x(t,\infty)$, which is a feasible  parallel parking motion.}

The AGHF being of parabolic type, it requires an {\it initial} condition (IC) $v(t) \in \X$
\begin{equation}\label{ic}
    x(t,0)=v(t),\quad t\in[0,T],
\end{equation}
as well as boundary conditions (BC), denoting the $i$-th state by $x_i$ \cite{WROCO2019}:
\begin{align}\label{bc}
    \mbox{if }x_i(0,s) \mbox{ fixed}&: x_i(0,s) = x_{\mathrm{init},i}\\
    \mbox{if } x_i(0,s) \mbox{ free}&: \dot x_i(0,s) = F_{d,i}(x(0,s)) \nonumber
\end{align}
The same holds for final conditions (at time $T$) with $\xfin$.

We denote by $x^\star(t)=x(t,\infty)$ the solution to which the AGHF converges. One can then extract controls from this curve that drive the system to the desired configuration.

To construct $G$, first, we pick a {\it bounded} $x$-dependent matrix $F_c(x)$, differentiable in $x$, so that the matrix
    \begin{equation}\label{def:bar_F}
    \bar F(x):=\begin{bmatrix}F_c(x)|F(x)
    \end{bmatrix} \in \R^{n \times n}
    \end{equation}
is invertible for all $x$. The matrix $F_c(x)$ can be found, e.g., via the Gram-Schmidt procedure. The columns of $F_c(x)$ are infinitesimal directions of motion that are not directly actuated. We then define the Riemannian metric tensor:
\begin{equation}\label{eq:G}
G(x):=(\bar F(x)^{-1})^\top D\bar F(x)^{-1}
\end{equation}
where $D:=\diag(\underbrace{\lambda,\cdots,\lambda}_{n-m},\underbrace{1,\cdots,1}_{m})$ for some large constant $\lambda>0$.
The parameter $\lambda$ can be thought of as a penalty on the infinitesimal  directions $F_c(x)$. Using the metric~\eqref{eq:G}, we can measure the ``actuated length'' of a curve $x(t)$ by integrating  $L(x,\dot x)$ as indicated earlier. But note that since $\lambda$ penalizes $F_c(x)$, a curve of minimal actuated length will use these directions only minimally, and this will yield a trajectory that the system can follow, with very high precision (quantitative relations between $\lambda$ and precision are derived in~\cite{AGHF2019}). {Note that the motion planning problem has to be feasible and the curve will converge to a local minimum, which is a feasible solution  in some neighbourhood of the initial condition.}

\section{Robot Dynamics} \label{dynamics}

\subsection{Single Rigid Body Model}\label{single_rigid_body_model}

We now apply these ideas to plan the motion of a 2D legged robot with massless legs--the Single Rigid Body Model, see Fig.~\ref{fig:massless_leg_model}. The torso of robot is a rigid body with mass $M$ and inertia $I$ around its CoM. The CoM position is $p=[p_x,p_y] \in \R^2$ and the torso's orientation is $\theta \in \mathbb{T}^1$. 

The robot has $k>0$ legs, each leg has point foot at the distal end, with coordinate $p_{i} = [p_{ix}, p_{iy}]^\top \in \R^2$. The contact force applied on the point foot at $p_i$ is $f_{i} = [f_{ix}, f_{iy}]^\top \in \R^2$.
The number of leg links and their lengths are not predetermined.  We  ensure the joint angles are feasible by adding constraints on the foot and hip positions; for example, if a leg has two links from hip to foot and one joint at knee, then the joint angle is feasible if the distance between foot and hip is smaller than the sum of link lengths. This kinematic constraint is discussed in Section \ref{constraints}.

The joint torques and contact forces can be mapped to each other by the foot Jacobians for the massless leg robot. Therefore we directly use the contact forces as inputs to the system. % and omit the mapping to joint torques. 
The equations of motion for the robot thus are:
\begin{equation}\label{eq:EoM}
 \ddot{p} = \frac{1}{M}\sum_{i=1}^k f_{i}-\begin{bmatrix}0 & g \end{bmatrix}^\top, \;
 \ddot{\theta} = \frac{1}{I}\sum_{i=1}^k f_i \times (p-p_i)
 \end{equation}
where $g$ is the gravitational acceleration.
For a leg $i$ with $\bar k_i$ total joints, when the foot is in contact, the joint torques $\tau_i = [\tau_1, \tau_2, ... \tau_{\bar k_i}]^\top$ can be easily mapped to contact force by using Jacobian  of the foot:
\begin{equation}\label{eq:force_torque}
    \tau_i = {J_i}^\top j_i
\end{equation} 
where $J_i \in \R^{2 \times \bar k_i}$ is the Jacobian of the $i$-th leg's foot to ground inertia frame. 
The terrain is {\bf not} assumed to be flat and is given as the zero-set of a $C^2$ function
\begin{equation}\label{eq:terrain}
     f_{\mathrm{terr}}(c_x,c_y) = 0
\end{equation}
where $[c_x, c_y]^\top$ is a point on the 2D terrain.%, and $f_{\mathrm{terr}}(c_x,c_y)$ is a $C^2$ terrain function on $[c_x, c_y]^\top$.

\subsection{Constraints for Legged locomotion}\label{constraints}
 The robot dynamics \eqref{eq:EoM} describes a rigid body controlled by external forces $f_i$ provided at contact points $p_i$. There are various constraints on $f_i$ and $p_i$ depending on whether there is contact with the ground. 
For each  leg, there are two different modes or phases: {1.~\bf Stance phase:} the foot is in no-slip contact with a surface and  {2.~\bf Flight phase:} the foot is in the air. We now describe the constraints in different phases. In \textbf{stance phase} for leg $i$:
%\begin{subequations}
%  \begin{tabularx}{\columnwidth}{p{3.2cm}X p{5.8cm}X}
%   \renewcommand{\arraystretch}{.8}
%  \begin{equation}
%    \label{eq:constr_contact}
%      \dot{p}_i = 0
%  \end{equation}
%  &\vspace{-.55cm}
%  \begin{equation}
%    \label{eq:constr_height}
%   f_{\mathrm{terr}}(p_{ix},p_{iy}) = 0
%  \end{equation}
%  \\ \vspace{-1.2cm}
%  \begin{equation}
%    \label{eq:constr_push}
%      f_i \cdot \overrightarrow{N}(p_{i}) \geq 0
%  \end{equation}
%  &\vspace{-1.2cm}
%  \begin{equation}
%    \label{eq:constr_cone}
%    |f_i \cdot \overrightarrow{T}(p_{i})| \leq \mu f_i \overrightarrow{N}(p_{i})
%\vspace{-.15cm}
%  \end{equation}
%  \end{tabularx}
%\end{subequations}
 \begin{itemize}
 \setlength\itemsep{-1em}
     \item[] \begin{equation}\label{eq:constr_contact}
       \dot{p}_i = 0
       \end{equation}
     \item[] \begin{equation}\label{eq:constr_height}
       f_{\mathrm{terr}}(p_{ix},p_{iy}) = 0
       \end{equation}
     \item[] \begin{equation}\label{eq:constr_push}
       f_i \cdot \overrightarrow{N}(p_{i}) \geq 0
       \end{equation}
     \item[] \begin{equation}\label{eq:constr_cone}
       |f_i \cdot \overrightarrow{T}(p_{i})| \leq \mu f_i \overrightarrow{N}(p_{i})
       \end{equation}
 \end{itemize}
the foot is in contact with the ground and has zero velocity:~\eqref{eq:constr_contact}-\eqref{eq:constr_height}.
The force $f_i$ is generated through contact with the ground, thus subject to the following constraints: let $\overrightarrow{N}(p_{i})$ and $\overrightarrow{T}(p_{i})$ be the unit normal and tangent vectors at the contact point $p_i$,  which can be directly calculated using the gradient of terrain function $f_{\mathrm{terr}}(\cdot,\cdot)$. The foot can only push against the ground: the projection of $f_i$ in normal direction to the surface at the contact point has to be positive~\eqref{eq:constr_push}. Furthermore, the contact force is constrained by a friction cone, formed by the normal contact force and the friction coefficient $\mu$, as shown by the light blue triangle in Fig.~\ref{fig:massless_leg_model}; see~\eqref{eq:constr_cone}.
In \textbf{flight phase} for leg $i$, there is no contact force and the feet are above the ground:
\begin{subequations}
  \begin{tabularx}{\columnwidth}{p{5.5cm}X p{5.8cm}X}
 \begin{equation}\label{eq:constr_flight}
      |f_i| = 0
      \end{equation}
  &
 \begin{equation}\label{eq:constr_foot_height}
       f_{\mathrm{terr}}(p_{ix},p_{iy}) \geq 0
      \end{equation}
  \end{tabularx}
  \end{subequations}
%\begin{equation}\label{eq:constr_flight}
%      |f_i| = 0
%      \end{equation}
% The contact force is zero during flight phase, with the absence of contact, e.g., $f_2$ in Fig.~\ref{fig:massless_leg_model}.
Finally, the following constraints enforce that the joint angles are feasible and hold for \textbf{both phases}:

\begin{subequations}
  \begin{tabularx}{.95\columnwidth}{p{5.1cm}X p{1.9cm}X}
 \begin{equation}\label{eq:constr_kinematics}
     |p-p_i| \leq R
     \end{equation}
  &
 \begin{equation}\label{eq:constr_py}
    f_{\mathrm{terr}}(p_x,p_y) \geq h_c
  \end{equation}
  \end{tabularx}
  \end{subequations}
  
Indeed, assuming that the hip for all legs are at CoM and the leg links are connected in series, the feasibility of joint angles can be ensured~\eqref{eq:constr_kinematics}, where $R$ is chosen to be less than the sum of the leg link lengths. One can make $R$ smaller to avoid approaching singularity configurations of the leg.If no collision with torso is desired, this constraint can be replaced by $(|p-p_i|-R)^2 \leq \Delta R^2$ with proper $\Delta R$. 
Constraint \eqref{eq:constr_py} ensures that the CoM is higher than the terrain height by some positive constant $h_c$.%{\color{blue}, which can be determined  so that the torso is always above the terrain}. %$h_c$ can be determined by the size of the torso.

\subsection{State Space Model}\label{state_space}

To represent the system in form of~\eqref{eq:sys}, define the state: 
$$x=[p, \theta, \dot{p}, \dot{\theta}, f_1, p_1, f_2, p_2,...f_k, p_k]^\top \in \R^{6+4k},$$
 in which the original controls $f_i$ and $p_i$ are states of the system, and  introduce the new controls: $ [u_i, v_i]^\top \!=\! [\dot{f}_i, \dot{p}_i]^\top \in \R^4$, which are the rate of change for the original controls $f_i$ and $p_i$. Now, the control to the system is $u = [u_1, v_1, \dots, u_k, v_k]^\top$. This operation allows us to encode constraints on the original controls as state constraints, as discussed in Sec.~\ref{constraints_state}. Denoting $m \times n$ zero matrix by $O_{m \times n}$ and $k \times k$ identity matrix by $I_k$, we can write the system in form of \eqref{eq:sys} with:
\begin{equation}\label{eq:state_space}
 F_d(x) = \begin{bmatrix} \begin{bmatrix}x_4&x_5&x_6\end{bmatrix}^\top \\ \frac{1}{M}\sum_{i=1}^k f_i - \begin{bmatrix}0 & g \end{bmatrix}^\top \\  \frac{1}{I}\sum_{i=1}^k f_i \times (p-p_i) \\ O_{4k\times 1} \end{bmatrix}, \; F(x) = \begin{bmatrix} O_{6 \times 4k} \\ I_{4k} \end{bmatrix}
 \end{equation}
where, from the definition of $x$ above,  $f_i = [x_{3+4i},x_{4+4i}]^\top$, $p_i = [x_{5+4i},x_{6+4i}]^\top$ and $p=[x_1,x_2]$. The 2-D cross product ``$\times$'' for torque calculation is defined as $[x_1, y_1]^\top \times [x_2, y_2]^\top = x_1 y_2 - x_2 y_1$. The drift term $F_d(x)$ includes all the robot dynamics, and the columns of $F(x)$ are the actuated directions, in other words, the directions that can be directly controlled by $u$.

%According to \eqref{def:bar_F}, we take $$F_c := [I_{6}, O_{4k \times 6}]^\top$$ so that $span\{F_c\}$ is orthogonal to $span\{F\}$, and $\bar{F} = I_{4k+6}$ is full rank for all $x$.

\section{A Riemannian metric for Legged Locomotion}\label{forlumation}

The motion planning problem is to find a trajectory for $x$ obeying~\eqref{eq:state_space} under constraints \eqref{eq:constr_contact} to \eqref{eq:constr_py}, for given boundary conditions, time span $T$.
 Up to this point, we have modeled the legged locomotion dynamics using continuous representations, that is, system \eqref{eq:sys} and constraints \eqref{eq:constr_contact} to \eqref{eq:constr_py}, in which \eqref{eq:constr_height} to \eqref{eq:constr_py} are state constraints  and \eqref{eq:constr_contact} is a constraint on input $\dot{p}_i = v_i$.
The challenge is that the constraints are different during flight and stance phases. In this work, the active phases are predefined, i.e.,  the timing of taking off and landing of each leg are set in advance. This results in a time-varying Riemannian metric. The phases can be also determined by the states, if the contact sequence is unscheduled, which would result in a non-smooth, but time-invariant, Riemannian metric. We will address this case in an upcoming publication.
The transition of a constraint from active to inactive is a discrete event, which we formulate  with the help of an {\it Activation Function} in Sec.~\ref{activation}. 
Then the input and state constraints are equipped with proper activation functions, which are discussed in Sec.~ \ref{constraint_fixed_contact} and Sec.~ \ref{constraints_state}.

\paragraph{\textit{\textbf{Activation Function}}} \label{activation}

\begin{figure}
 \tikzset{every picture/.style={scale=1}}
   \centering
   \scalebox{0.6}{
   \includegraphics{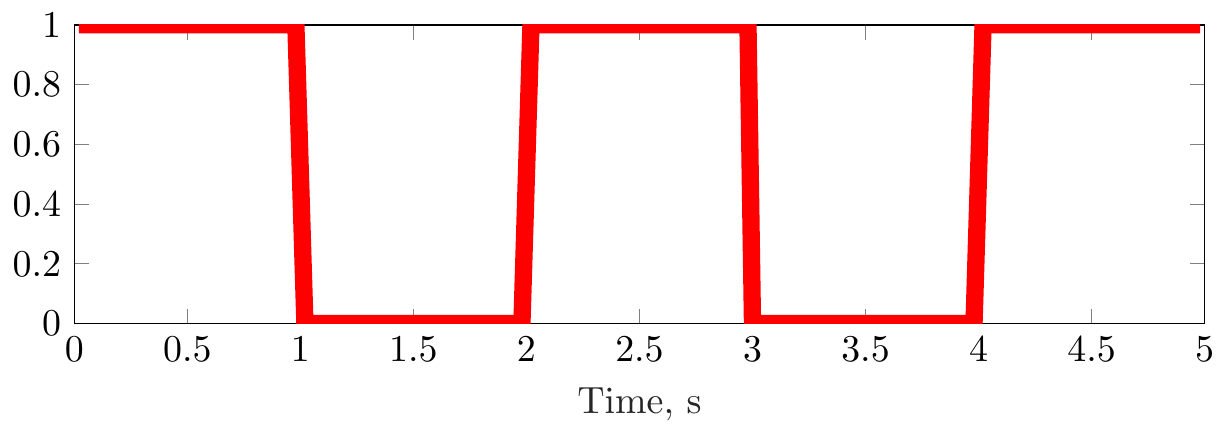}
   }
 \caption{\small A switch function for $K=3$, the stance phases are $t=[0,1]\cup [2,3]\cup [4,5]$.}
 \label{fig:switch_function}
 \end{figure}

The {\it Activation Function} $A_i(t,x)$ for the $i$th leg is  is a binary valued function which determines whether the leg is in stance or flight:
\begin{equation}\label{eq:activation}
A_i(t):= \begin{cases}
    1, & \text{if foot in stance}\\
    0, & \text{if foot in flight}.
  \end{cases}
\end{equation}
We are given the time sequences $\{t_{i,1}^j\}_j$ and $\{t_{i,2}^j\}_j$ of landing and take off time of step $j$ for leg $i$ respectively. The activation function $A_i(t)$ is defined as:
\begin{equation}\label{eq:activation_pre}
A_i(t):= \sum_{j=1}^k H(t-t^j_{i,1})-H(t-t^j_{i,2}).
\end{equation}
where $k$ is the total number of steps of the foot and $H(c)$ is a Heaviside unit step function. As a result, if the foot is in stance phase, the value of switch function is 1, otherwise 0. An example of switch function is shown in Fig.~\ref{fig:switch_function}.

\paragraph{\textit{\textbf{Fixed Contact Foot Position Formulation}}} \label{constraint_fixed_contact}
 The position of the contact point of the foot with the ground is constant during stance phase~\eqref{eq:constr_contact}. This constraint is an input constraint since the foot velocity $\dot p_i = v_i$ is an input to the state space model~\eqref{eq:sys}. However, $v_i$ is either free or zero depending on the phase of the foot. 
To encode the switching between free and constant foot velocity, i.e. flight and stance phases, we  define the Riemannian metric $G$ similarly to~\eqref{eq:G}, with a {\it time varying} penalty matrix $D(t)$ equipped with activation functions of all legs:
\begin{equation}\label{eq:D_fixed}
D=\diag(\underbrace{\lambda,\cdots,\lambda}_{6},\Lambda_1, \cdots,\Lambda_k)
\end{equation}
with $\Lambda_i\!=\!\diag \big(1,1,1\!+\!\lambda A_i(t),1\!+\!\lambda A_i(t)\big)$, where $A_i(t)$ is the activation function for leg $i$. The  term $1\!+\!\lambda A_i(t)$ in $\Lambda_i$ penalizes the length of trajectories that use the control $v_i$. If the foot is in stance, the value of $A_i$ is 1 and the curve length is increased due to a nonzero $v_i$ multiplied $\lambda$; as a result, $v_i$ will be minimized.  When the foot is in flight phase, $A_i$ is 0, therefore the value of $v_i$ is free as it won't affect the curve length. Hence curves of minimal length for this metric are so that  the constraint~\eqref{eq:constr_contact} is met when the foot is in stance.

With the  penalty terms for $u_i$ set to 1,  we are not constraining the changing {\it rate} of contact forces, $\dot{f}_i$, in the Riemannian metric. The constraints on $f_i$ (cone~\eqref{eq:constr_cone} and positivity~\eqref{eq:constr_push} constraints)  will be encoded in Sec.~\ref{constraints_state} as state constraints.  % The force constraints are more complicated due to the constraints on its magnitude and direction, \eqref{eq:constr_cone} and \eqref{eq:constr_push}, which can be better handled by states constraints formulation.
By the formulation \eqref{eq:D_fixed}, the control is either free or constrained to zero depending on the time dependent penalty matrix $D(t)$. However, one can also have constraints on the magnitude of control by letting $u$ be a state of the system which is directly controlled by a newly introduced unconstrained control. The constraint on $u$ is then a state constraint. This is exactly the intention of formulating contact forces and foot positions as states in Sec. \ref{state_space}.

Owing to the simple structure of $F$ in~\eqref{eq:state_space}, we take $F_c \!:=\! [I_{6}, O_{4k \times 6}]^\top$ so that $\operatorname{span}\{F_c\}\!\perp\! \operatorname{span}\{F\}$ and $\bar{F} \!=\! I_{4k+6}$ is full rank for all $x$. 

\paragraph{\textit{\textbf{State Constraints Formulation}}} \label{constraints_state}
Each of the state constraints from \eqref{eq:constr_height} to \eqref{eq:constr_py} can be formulated as scalar function $h(x) \!=\! 0$ for equality constraints or $h(x) \!\leq\! 0$ for inequality constraints. For example, for constraint \eqref{eq:constr_py} we have $h\!=\!h_c\!-\!f_{\mathrm{terr}}(p_x,p_y)\!\leq\! 0$. We modify slightly the method in~\cite{IFAC_Mechatronics2019} to encode the activation/deactivation of the state constraints in different phases. To this end,  denote by  $h_j(x)$, $j \!\in\! \z^+ \!\leq\! k_c$, the $j$-th scalar constraint function with $k_c$  the  number of such constraints. We apply the following method to {\it each} constraint individually.

First we add one state  per constraint, denoted by $\zeta_j$, resulting in the {\bf augmented state} of the system: $\hat{x}\!:=\![x^\top, \zeta_j]^\top$. The new state $\zeta_j$ keeps track of the {\it accumulated signed}  error between $h_j(x)$ and zero: $\zeta_j(t)\!:=\!\int_0^t h_j(x(\tau))S_j(\tau,x(\tau)) d\tau \!\longrightarrow \!\dot{\zeta}_j\!:=\! h_j(x(t))S_j(t,x(t))$, where $S_j(t,x)$ is a scalar {\it switch function} for different type of constraints:
\begin{equation}\label{eq:switch_function}
S_j(t,x) = \begin{cases}
    B_j(t), &\!\!\! \text{for equality constraint}.\\
    H(h_j(x))B_j(t), &\!\!\! \text{for inequality constraint}.
  \end{cases}\!\!\!\!\!\!
\end{equation}
in which $H(\cdot)$ is the  Heaviside  function used in~\eqref{eq:activation_pre},  $B_j(t)$ is the {\it activation function} for the $j$-th scalar state constraint. A constraint can be active either for the duration of the motion (e.g., bound on joint angle), or depending on the phase in the motion. The activation function for the $j$-th scalar state constraint is thus given by:
\begin{equation}\label{eq:activation_function_const}
B_j(t)\!=\!\begin{cases}
    1 & \!\!\!\!\text{if holds all the time}\\
    A_i(t) & \!\!\!\!\text{if holds in stance phase of foot $i$}\\
    1 \!- \!A_i(t) & \!\!\!\!\text{if holds in flight phase of foot $i$}.
  \end{cases}
\end{equation}
For example, the activation function for constraint \eqref{eq:constr_flight} of the $i$-th leg is $1-A_i(t)$.

It is now easy to see that thanks to the switch function~\eqref{eq:switch_function}, the new state $\zeta_j$ is zero if the constraint is satisfied during $t=0$ to $t=T$.
 The augmented system dynamics is thus given by
\begin{equation}\label{eq:state_space_const}
\dot{\hat{x}} = \hat{F}_d(\hat{x}) + \hat{F}(\hat{x})\begin{bmatrix} u \\ h_j(x)S_j(t,x)\end{bmatrix}
\end{equation}
where
\begin{equation}
\hat{F}_d(\hat{x}) =  \begin{bmatrix} F_d(x) \\ 0  \end{bmatrix}, \;
\hat{F}(\hat{x}) =  \begin{bmatrix} F(x) & O_{6+4k,1} \\ O_{1,6+4k} & 1 \end{bmatrix},
\end{equation}
in which $\hat{F}_d$ is the augmented drift.

\paragraph{\textit{\textbf{Actuated Curve Length and Geometric Heat Flow}}}\label{Actuated_and GHF}

The Riemannian metric $G(t)$ defined by \eqref{eq:G} and \eqref{eq:D_fixed} is constructed so that the length of paths violating the constraints is large. We do so by introducing
\begin{equation}\label{eq:G_hat}
\hat{G}(t):=\begin{bmatrix} G(t) & O_{6+4k,1} \\ O_{1,6+4k} & \lambda_j  \end{bmatrix}
\end{equation}
where $G(t)$ is the Riemannian metric defined for the original states $x$ in \eqref{eq:G} and $\lambda_j$ is a large constant penalizing the infinitesimal directions that violate the constraints. With this metric, we obtain that the actuated curve length of the augmented system is
\begin{align}\label{eq:actuated_length_new}
\hat{L}&=\int_0^T \big({(\dot{\hat{x}}-\hat{F}_d(\hat{x}))^\top \hat{G}(t) (\dot{\hat{x}}-\hat{F}_d(\hat{x}))}\big)^{1/2} dt 
\end{align}
\begin{comment}
For the augmented system, we define the {\it constrained AGHF}:
\begin{equation}\label{eq:aghf_const}
    \frac{\partial x(t,s)}{\partial s} = \nabla_{\dot x(t,s)}\left(\dot x(t,s)-F_d\right)+r(t,s)-\lambda_j \frac{\partial (h(x)^2 S(t,x)) }{\partial x}
\end{equation}
where the first two terms of the right-hand side are the same as~\eqref{eq:aghf}.
%The added term $-\lambda_j \frac{\partial (h(x)^2 S(x))}{\partial x}$ acts on the curve to reduce the length of $\zeta$.+++++ This term involves taking derivative of the switch function, therefore the smooth approximation of step function is used in \eqref{eq:switch_function}.
Similarly as in Section \ref{AGHF}, the constrained AGHF \eqref{eq:aghf_const} minimize the actuated curve length of the augmented system: 

\begin{align}\label{eq:actuated_length_new}
&\int_0^T \big({(\dot{\hat{x}}-\hat{F}_d(\hat{x}))^\top \hat{G}(t,\hat{x}) (\dot{\hat{x}}-\hat{F}_d(\hat{x}))}\big)^{1/2} dt \nonumber \\ 
&= \int_0^T\big( (\dot x-F_d(x))^\top G(t,x) (\dot x-F_d(x)) + \lambda_j h(x)^2 S(t,x) \big)^{1/2} dt
\end{align}
\end{comment}
and plugging in the defining of $G$, we get $$\hat L=\int_0^T
\big((\dot x-F_d(x))^\top G(t) (\dot x-F_d(x)) 
+ \lambda_j h_j^2(x) S_j(t,x) \big)^{1/2} dt.$$ Compared to~\eqref{actuated_length}, it contains the additional term $\lambda_j h_j(x)^2 S_j(x)$. When the constraint $h_j(x)$ is active, namely, $B_j(t)=1$: for equality constraint $h_j(x)=0$, by construction $S_j(x)=1$, the actuated length is penalized by $\lambda_j$ when $h_j(x)$ is not close to zero, driving $h_j(x)\rightarrow 0$. For an inequality constraint $h_j(x)\leq 0$, the actuated length is penalized by $\lambda_j$ when $h_j(x) > 0$ and $S_j(x)=H(h_j(x))\approx 1$, driving $h_j(x)\rightarrow 0$. However, this term has no effect on the actuated length if $h_j(x) \leq 0$ since $S_j(x)=H(h_j(x))\approx 0$. Hence,  minimizing the actuated length of an augmented state trajectory results in {\it minimizing the violation} of the state constraints while it is active. 
In conclusion, solving~\eqref{eq:AGHF_general} derived from Lagrangian~\eqref{eq:actuated_length_new} leads to a curve with minimum actuated length, which is a motion admissible for system~\eqref{eq:EoM} and respects the constraints~\eqref{eq:constr_contact}-\eqref{eq:constr_py}.

\paragraph{\textit{\textbf{Step Function Approximation}}}\label{step_approx}
Since the AGHF~\eqref{eq:AGHF_general} requires to take  derivatives of a Heaviside step function (whose derivative formally does not exist as a function),  we approximate the step function $H(c)$ and its derivative $\frac{dH}{dc}(c)$ using  a logistic approximation~\eqref{eq:step_approx} (below).
%\begin{equation}\label{eq:step_approx}
%H(c):= \frac{1}{1+e^{-\alpha c}}
%\end{equation}
It approximates a unit step at $c=0$, namely, $H(c)\approx0$ if $c<0$ and $H(c)\approx1$ if $c>0$. The constant $\alpha$ controls the accuracy of the approximation. The derivative is $\frac{\alpha e^{-\alpha c}}{(1+e^{-\alpha c})^2}$, which  causes overflow problem when evaluating its value numerically for large $\alpha$. To circumvent this numerical issue, we use zero-centered normal distribution~\eqref{eq:step_dot_approx} to approximate the derivative of step function,
%\begin{equation}\label{eq:step_dot_approx}
%\frac{dH}{dc}:= \frac{1}{\beta \sqrt{2 \pi}} %e^{-(c/\beta)^2}
%\end{equation}
where $\beta$ is a large number scaling the value of $\dot{H}(c)$ at $c=0$.

\begin{subequations}
  \begin{tabularx}{.95\columnwidth}{ 
   >{\raggedright\arraybackslash}X 
   >{\raggedright\arraybackslash}X  }
  %{p{3.2cm}{}X p{7cm}X}
%   \renewcommand{\arraystretch}{.8}
 \begin{equation}\label{eq:step_approx}
\hspace*{-1.6cm} H(c):= \frac{1}{1+e^{-\alpha c}}
\end{equation}
 &
\begin{equation}\label{eq:step_dot_approx}
\hspace*{-.28cm}\frac{dH}{dc}:= \frac{1}{\beta \sqrt{2 \pi}} e^{-(c/\beta)^2}
\end{equation}
\end{tabularx}
\end{subequations}

\paragraph{\textit{\textbf{Control Extraction}}}\label{control_extraction}
Having solved numerically the  AGHF~\eqref{eq:AGHF_general} with BC~\eqref{bc} and IC~\eqref{ic}, we denote the steady state solution by $x^*(t) \!:=\! x(t,s_{max})$ for large enough $s_{max}$.
The control can be extracted from it according to:
\begin{equation}\label{eq:control_extraction1}
u(t):= \begin{bmatrix}
O_{4k,6} & I_{4k} \end{bmatrix}\bar F^{-1}(\dot x^*(t)-F_d(x^*(t))
\end{equation}
This extracted control  drives the system {\it arbitrarily close} to the desired final state. When integrating~\eqref{eq:sys} with control \eqref{eq:control_extraction1}, a trajectory $\tilde x(t)$ is obtained, which is our solution to the trajectory planning problem. We call it {\bf integrated path}.

%This extracted control will drive the system arbitrarily close to the desired final state, {\it provided} that  a solution to the motion planning problem exists (which is, in general, an open problem for control systems with drift). When integrating~\eqref{eq:sys} with control \eqref{eq:control_extraction1}, a trajectory $\tilde x(t)$ is obtained, which is our solution to the trajectory planning problem. We call it {\bf integrated path}.

\paragraph{\textit{\textbf{Planning Algorithm Summary}}}\label{algorithm_summary}

\begin{comment}
\begin{figure*}
%\tikzset{every picture/.style={scale=1}}
  \centering
  %\input{figures/diver.tex}
  \scalebox{0.8}{\input{figures/flow_chart.tex}}
\caption{\small Illustration of the GHF planning algorithm architecture for legged locomotion}
\label{fig:flow}
\end{figure*}
\end{comment}

The steps of the algorithm can be summarized as follow:
\begin{description}
    \item[Step 1:] Specify number of legs $k$, find $F_d(x)$ and $\bar{F}$.
    \item[Step 2:] Specify timing for stance/flight phases of each leg, construct all $A_i(t)$ to build $D(t)$. Then construct metric $G(t)$ with $\bar{F}$ and $D(t)$ from \eqref{eq:D_fixed}, using large $\lambda$.
    \item[Step 3:] Specify the terrain $f_{terr}(\cdot,\cdot)$, $\mu$, $R$, $h_c$. Construct switch functions $S_j(t)$ for the state constraint using $A_i(t), B_j(t)$. Formulate all state constraints $h_j(x)$ for \eqref{eq:constr_height}--\eqref{eq:constr_py} as extra states equipped with $S_j(t)$ to get the augmented states $\hat x$.
    \item[Step 4:] Get the augmented drift $\hat F_d$ from $F_d$, construct the metric $\hat G$ \eqref{eq:G_hat}  using $G$ and large $\lambda_j$.
    \item[Step 5:] Find the actuated length \eqref{eq:actuated_length_new} using $\hat G$, $\hat F_d$, then find the corresponding AGHF~\eqref{eq:AGHF_general}.
    \item[Step 6:] Solve AGHF with BC~\eqref{bc} and IC~\eqref{ic} and large enough $s_{max}$ to obtain the solution $x^*(t)$.
    \item[Step 7:] Extract control $u(t)$ from $x^*(t)$ using \eqref{eq:control_extraction1}. Integrate the dynamics \eqref{eq:state_space} with control $u(t)$ and initial value $x^*(0)$ to obtain the integrated path $\tilde{x}(t)$, which is the planned motion.
\end{description}
\section{Application Examples}\label{applications}
To illustrate our method, we use it to plan the two different motions: one leg hopping on a flat terrain, and two leg hopping on an uneven terrain.
\subsection{One Leg Hopping on Even Terrain}\label{one_leg}
\begin{figure}
 \centering
 \scalebox{.8}{
 \includegraphics{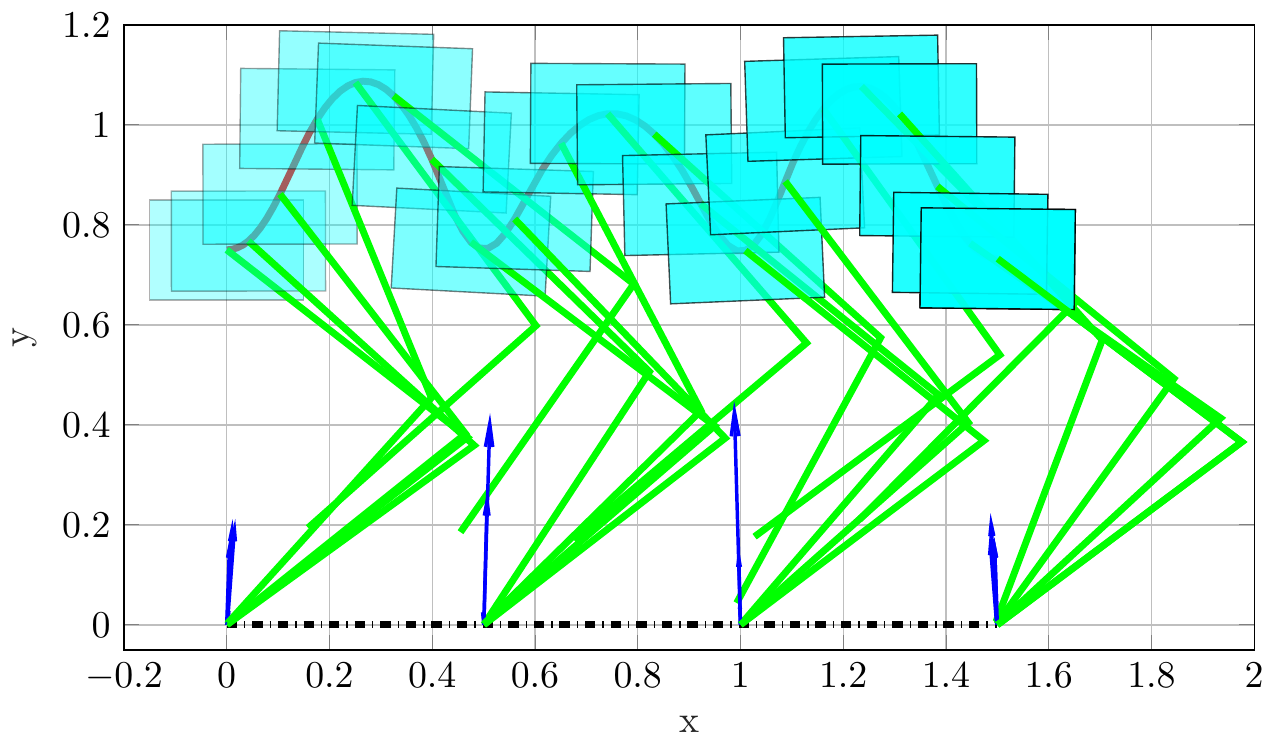}
 }
\caption{\small Snapshots of one leg hopping}
\label{fig:oneleg_snapshots}
\end{figure}

\begin{figure}
 \centering
 \scalebox{1}{
 \includegraphics{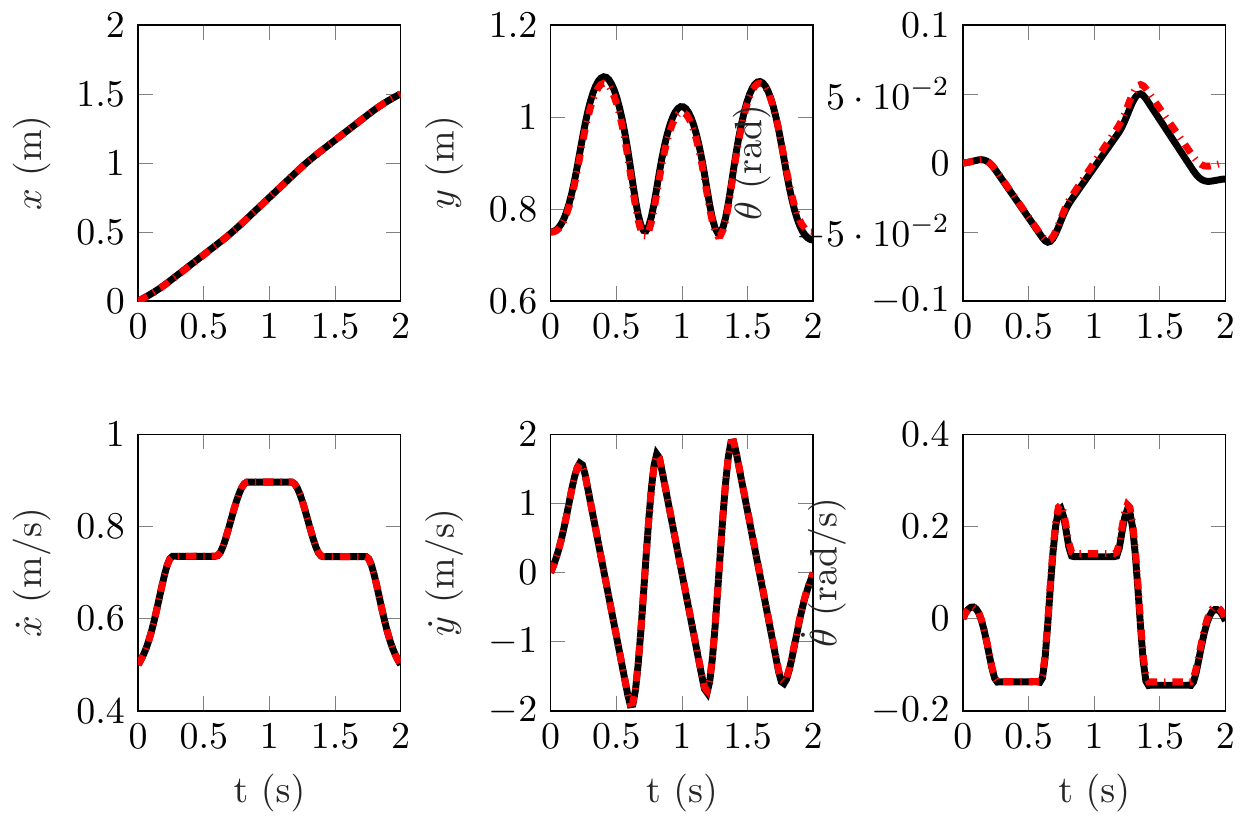}
 }
\caption{\small Torso position and orientation for one leg hopping}
\label{fig:states1}
\end{figure}

\begin{figure}
 \centering
 \scalebox{0.8}{
 \includegraphics{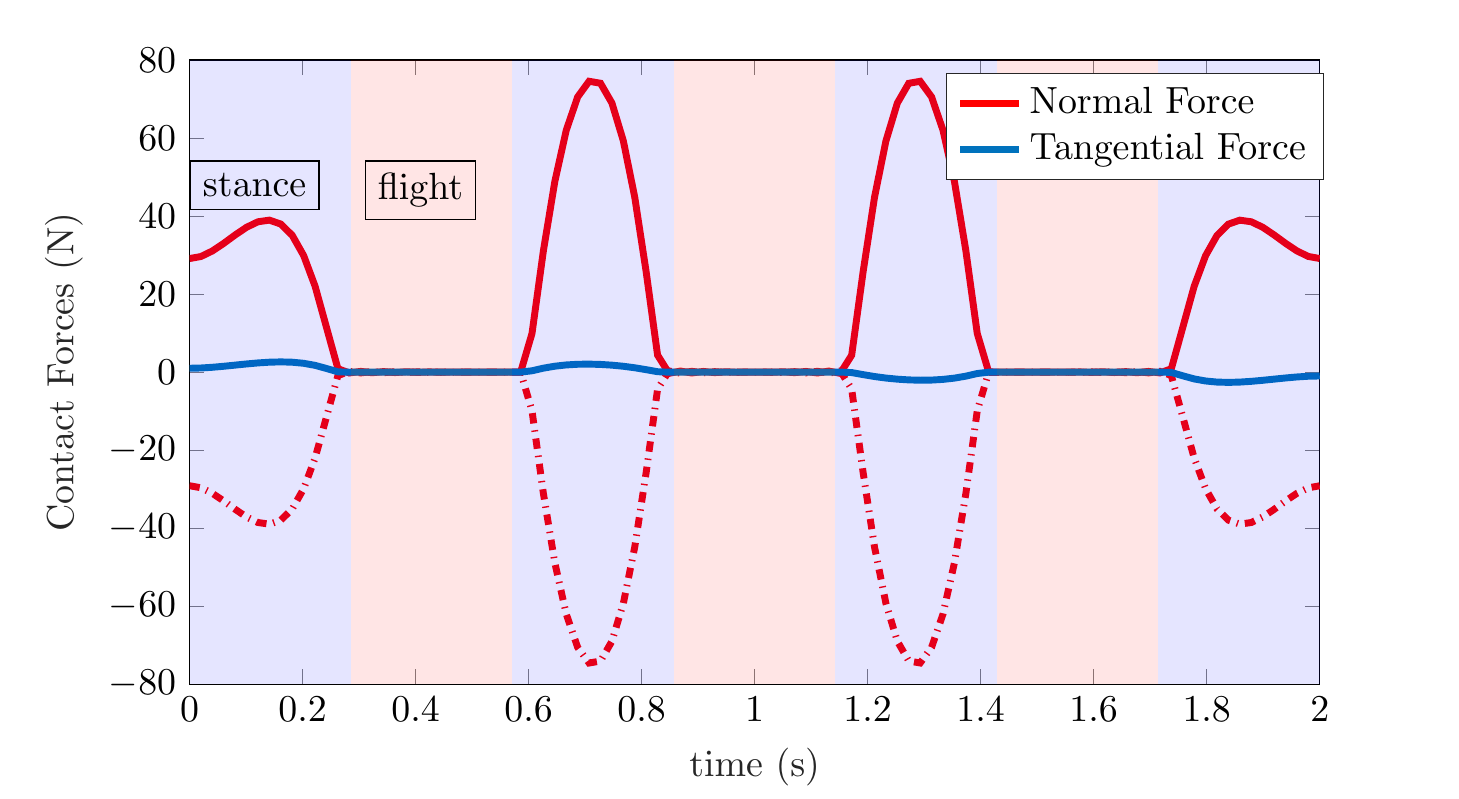}
 }
\caption{\small Contact force profile for one leg hopping}
\label{fig:oneleg_force}
\end{figure}

The first example to illustrate the method is one leg hopping motion. The terrain is flat therefore $f_{\mathrm{terr}}(c_x,c_y) = c_y$. The friction coefficient is chosen to be rubber to concrete friction coefficient, $\mu=1$. The leg is rooted at the CoM of torso. We choose the maximum radius $R$ in constraint \eqref{eq:constr_kinematics} to be $1$ $m$. The CoM is expected to be higher than the terrain by $0.3$ $m$ so $h_c=0.3$ in~\eqref{eq:constr_py}.
\begin{table}
\centering
 \begin{tabular}{||c | c ||} 
 \hline
 number of legs $k$ & $1$  \\
 \hline
 number of hops & $3$ \\
 \hline
 $T$ & $2$ $s$  \\
 \hline
 $\xin$ & $[0,0.75,0,0.5,0,0,\cdot,\cdot,0,0]^\top$  \\
 \hline
 $\xfin$ & $[1.5,0.75,0,0.5,0,0,\cdot,\cdot,1.5,0]^\top$  \\
 \hline
\end{tabular}
\caption{One leg hopping setup}
\label{tab:bc_oneleg}
\end{table}
The goal is to move the CoM from $[0,0.75]^\top$ to $[1.5,0.75]^\top$ with $T=2$ $s$ and $3$ hops. The stance/flight phases timing are predefined to have a $1:1$ ratio. Details of the setup of the motion is listed in Tab.~\ref{tab:bc_oneleg}, in which the symbol ``$\cdot$'' indicates the state that has free boundary value. The initial condition to solve the AGHF is a straight line  connecting  $\xin$ and $\xfin$ in time span $T$, which is not a feasible trajectory. Choosing the penalty $\lambda=2\times10^6$ and solving the AGHF \eqref{eq:AGHF_general} with a large $s_{max}$, we can have a steady state solution $x^\star (t)$. Using the extracted control \eqref{eq:control_extraction1} to integrate the system dynamics \eqref{eq:sys} with \eqref{eq:state_space}, the integrated path $\tilde x(t)$ can be obtained, which is the planned motion. Fig.~\ref{fig:oneleg_snapshots} shows snapshots of the motion at different time instants, and Fig.~\ref{fig:states1} shows the trajectories of torso position and orientation. In this latter figure, the red dotted line indicates the AGHF solution $x^\star$ and the black line indicates the integrated path $\tilde x$. We can observe that $x^\star(T)\approx\tilde x(T)$ with small deviation. In fact, the planning error $|x^\star(T)- \tilde x(T)|\rightarrow 0$ as $\lambda \rightarrow \infty$. The contact forces (normal and tangential) are shown in Fig.~\ref{fig:oneleg_force}; we observe that they are indeed zero during the flight phase. Also note that the tangential force is higher early in the motion, as it needs to impart the body with momentum, and negative at the end of the motion, as it needs to slow down the body.

\subsection{Two Leg hopping on uneven terrain}
\begin{table}
\centering
 \begin{tabular}{||c | c ||} 
 \hline
 number of legs $k$ & $2$  \\
 \hline
 number of hops & $3$ \\
 \hline
 $T$ & $2$ $s$  \\
 \hline
 $\xin$ & $[0,0.85,0,0.5,0,0,\cdot,\cdot,0,\cdot,\cdot,\cdot,0,\cdot]^\top$  \\
 \hline
 $\xfin$ & $[2,0.85,0,0.5,0,0,\cdot,\cdot,2,\cdot,\cdot,\cdot,2,\cdot]^\top$  \\
 \hline
\end{tabular}
\caption{Two leg hopping setup}
\label{tab:bc_twoleg}
\end{table}

\begin{figure}
\centering
\includegraphics[width=.75\columnwidth]{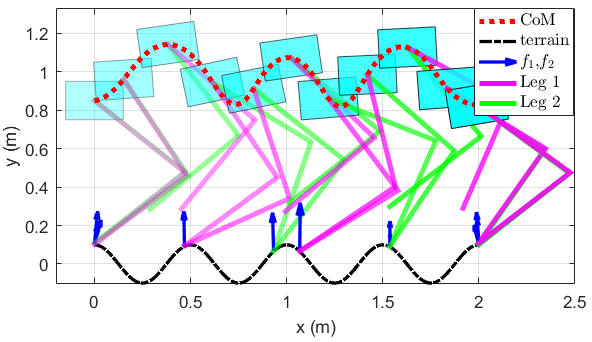}
\caption{\footnotesize Snapshots of two leg hopping on the terrain  $f_{\mathrm{terr}}(c_x,c_y)=c_y-0.1\cos(4\pi c_x)$. The torso (cyan box) has a mass $M\!=2\!$ $kg$ and inertia $I\!=\!1$ $kg\cdot m^2$. The leg is rooted at the CoM of torso and the maximum radius $R$ in constraint~ \eqref{eq:constr_kinematics} is $1$ $m$. The CoM is higher than $h_c=0.3$ in \eqref{eq:constr_py}. Friction coefficient is $\mu = 1$. The goal is to move the CoM from $[0,0.85]^\top$ to $[2,0.85]^\top$ with $T=2$ $s$ and $3$ hops. The stance/flight phases timing for leg 1 are predefined to have a $1:1$ ratio while the timing of the second leg is shifted by $-0.05$ $s$, so that the two legs are not synchronized. The generated motion shows the following gaits: Leg 1 kicks off at the first hump, then steps on the second and third hump while leg 2 kicks off at the first hump and then steps on the third and fourth hump. Both feet land on the last hump at the end.}
\label{fig:twoleg_snapshots}
\end{figure}

\begin{figure}
 \centering
 \scalebox{1}{
 \includegraphics{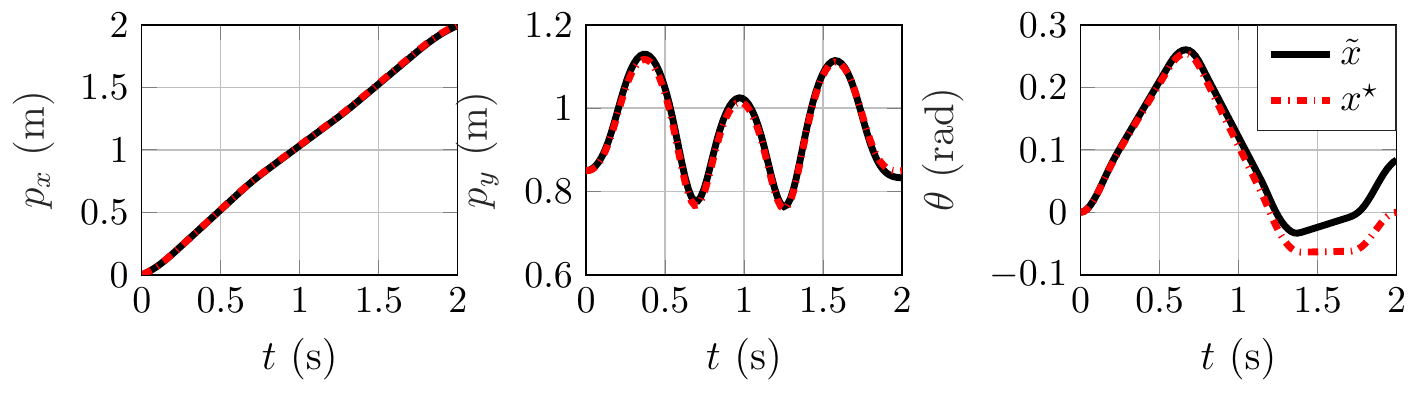}
 }
\caption{\footnotesize Torso position and orientation trajectories for AGHF solution $x^\star$ and integrated path $\tilde x$. }
\label{fig:states}
\end{figure}

We illustrate the performance of our method in planning motion for a two leg hopping robot  on uneven terrain with sinusoidal profile, see Fig.~\ref{fig:twoleg_snapshots} for snapshots of the planned motion\footnote{Supplementary animations can be found in playlist: \url{https://www.youtube.com/playlist?list=PLRi8ecX8Oy0UZzYoRi2DASQfU16sKtRph}}. The initial condition to solve the AGHF is a straight line  connecting  $\xin$ to $\xfin$ which is not a feasible trajectory. Choosing the penalty $\lambda\!=\!10^6$ and solving the AGHF \eqref{eq:AGHF_general} with a large $s_{max}\!=\!0.0005$ in \textrm{pdepe} from the Matlab PDE toolbox , we obtain $x^\star (t)$. Using the extracted control \eqref{eq:control_extraction1} to integrate the system dynamics \eqref{eq:sys} with \eqref{eq:state_space}, the planned motion $\tilde x(t)$ can be obtained. Fig.~\ref{fig:states} shows the trajectories of torso position and orientation. The trajectory is not necessarily periodic for each step since the timing of contact for each foot is not synchronized. The joint angles are guaranteed to be feasible since the distances between feet and the hip are constrained to be less than %{\color{blue}$1$ $m$ while }
the maximum distance each foot can reach. %{\color{blue} is $1.2$ $m$ (sum of leg links' length)}.
The algorithm is able to find a motion with different contact points for each leg while enforcing all the constraints. The contact points being automatically chosen, even in the case of uneven terrain, is quite advantageous compared to methods requiring  the user to predefine the contact positions or safe contact regions, as suboptimal choices here can strongly decrease the probability to find feasible solutions. However, this is done at the expense of a terrain function that needs to be $C^2$.

\begin{figure}
 \centering
 \scalebox{1}{
 \includegraphics{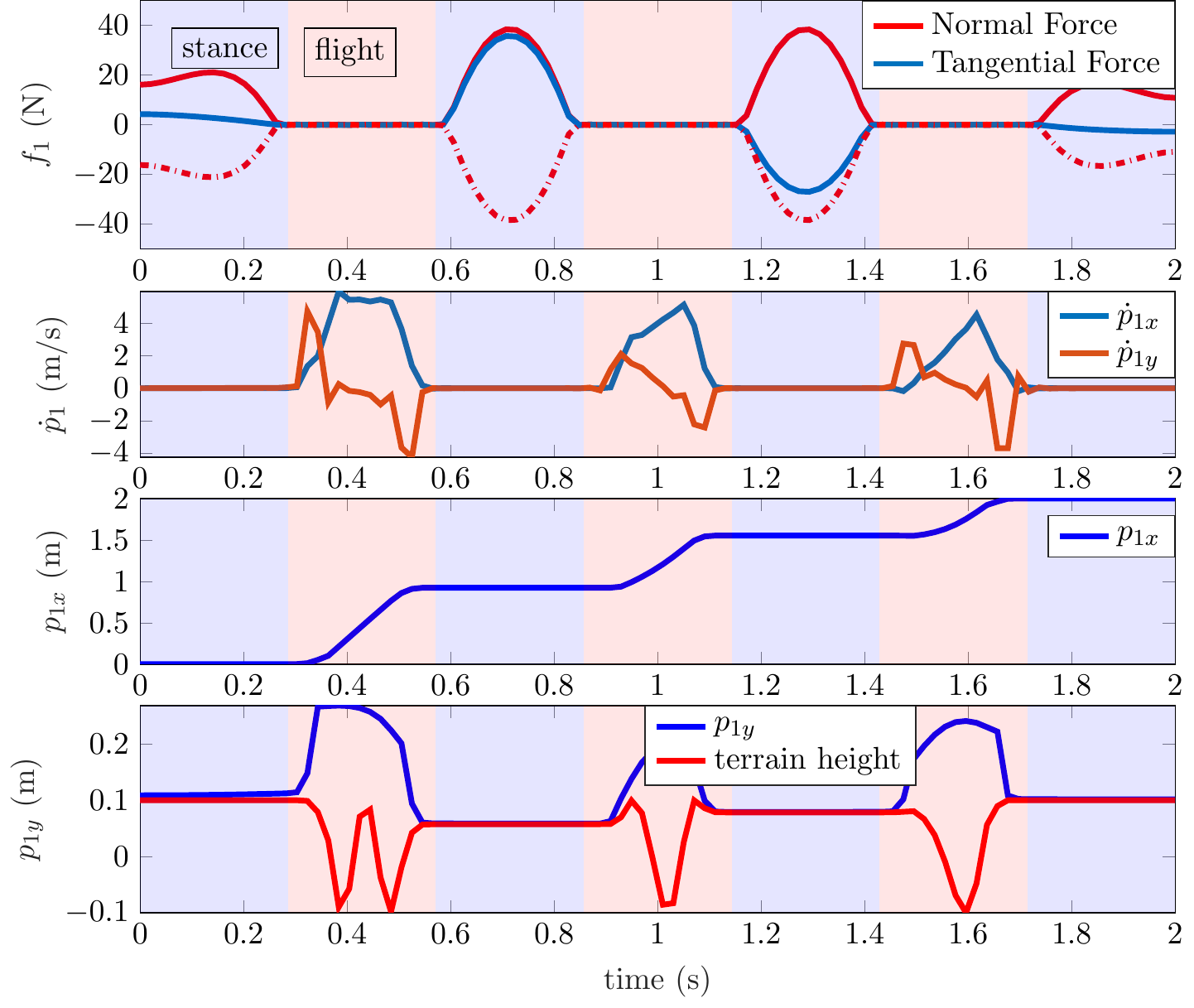}
 }
\caption{\small Foot contact force (top), velocity (second) and position (third and bottom) of leg 1}
\label{fig:twoleg_foot}
\end{figure}

What makes hopping on uneven terrain challenging is that the foot might slip on inclined surface. With the friction cone constraints enforced, we can see that feet are able to step on the inclined contact points where enough friction can be provided. With the feet stepping on inclined surface, the values of tangential forces are close to the boundary of the friction cone constraint, but are strictly inside of the friction cone boundaries enclosed by solid and dotted red lines.
Fig.~\ref{fig:twoleg_foot} top shows the contact forces for leg 1. The predefined stance/flight phases are indicated by light blue/red backgrounds. Both normal and tangential force being zero in flight phases reflects the constraint \eqref{eq:constr_flight}. In each stance phases, the normal force is larger than zero, as in~\eqref{eq:constr_push}. With $\mu\!=\!1$, The friction cone \eqref{eq:constr_cone} is  $|f_i \cdot \overrightarrow{T}(p_{i})| \leq f_i \overrightarrow{N}(p_{i})$. Namely, the positive normal force and its opposite value (red dotted line) serve as upper and lower bound for tangential force. It is clear in Fig.~\ref{fig:twoleg_foot} top, the value of tangential force is close to but strictly inside the boundary of the friction cone.

The fixed foot position constraint \eqref{eq:constr_contact} during stance phase is also achieved as shown in Fig.~\ref{fig:twoleg_foot}-2, where the foot velocity $\dot{p}_1$ is zero during stance phases. This results in constant foot position, as shown in Fig.~\ref{fig:twoleg_foot} third and bottom for $x$ and $y$ component of foot position $p_1$. Furthermore, in Fig.~\ref{fig:twoleg_foot} bottom, the foot height $p_{1y}$ is constrained to be the same as terrain height to make sure the foot is in contact with ground during stance. This reflects the constraint \eqref{eq:constr_height}. For flight phases, the algorithm is able to plan the smooth swing motion trajectory of the foot by actuating $p_1$ with proper $\dot{p}_1$.

\section{Discussion and Conclusion}\label{discussion_conclusion}

\begin{figure}
 \centering
 \scalebox{1}{
  \includegraphics{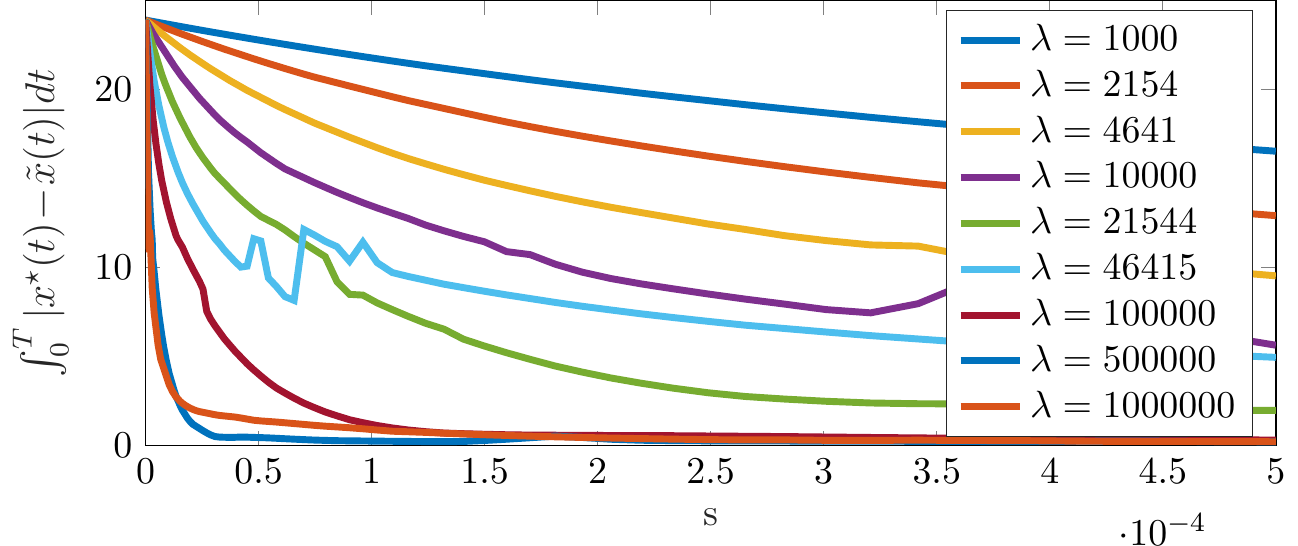}
}
\caption{\small Convergence of the planning error $e$ along $s$ and $\lambda$ for the 2 leg hopping example}
\label{fig:s_vs_error}
\end{figure}

\textbf{\textit{Parameter selection and convergence:}} We now elaborate on convergence properties and parameter selection.% for the algorithm.
%The selection of $\lambda$, $s_{max}$, $\alpha$ and $\beta$ are critical for finding a feasible solution. 
A large  $s_{max}$ is used to obtain the steady state AGHF solution, whereas a large $\lambda$ guarantees that error between AGHF solution and the integrated path, denoted by $e\!=\!\int_0^T |x^\star(t)\!-\!\tilde x(t)|dt$, is small. The planning error $e$ can serve as a measure of violation of dynamics and constraints, and it can be shown to converge to zero as $\lambda,\!s\!\rightarrow\!\infty$, see Fig.~\ref{fig:s_vs_error}. Given large enough $\lambda$, the steady state AGHF solution is a feasible trajectory for the approximated robot dynamics instead of the actual hybrid robot system. The model accuracy increases with larger $\alpha$ and $\beta$, resulting in a better approximation of the discrete transition between flight and stance phase. %That is, with sufficiently large $\alpha$ and $\beta$, the violation of constraints near the switching time is negligible.

\noindent \textbf{\textit{Computational efficiency:}}
The complexity of solving the AGHF numerically is low, especially compared to other PDEs in motion planning based on the Hamilton-Jacobi-Bellman (HJB) equation. The reason for this is that we {\it always} solve, regardless of the dimension $n$ of the system, for a function with a {\it domain} of dimension 2 (with variables $(t,s)$). 
Thus, whereas PDEs such as the HJB scale {\it exponentially} in the dimension of the system, the AGHF scales {\it polynomially}. %{\color{blue}(in the order of matrix multiplication)}. 
However in our implementation, choice of larger $\lambda$, $\alpha$, $\beta$ can result in longer solving time, since the \textrm{pdepe} function determines the integration step size $\Delta s$ automatically 	without option to customize \cite{IFAC_Mechatronics2019}. For example, with large value of $\lambda$, the magnitude of $\frac{dx(\cdot, s)}{ds}$ is generally large, as a result, the solver uses smaller $\Delta s$, increasing the solving time. However the heat flow being a parabolic PDE, it can be solved in parallel~\cite{horton1995algorithm}, thus potentially vastly improving the solving time.

\noindent \textbf{\textit{Conclusion and future work: }} The framework introduced in this paper allows us to encode motion planning for legged robot in a simple and unified way: a Riemannian inner product is defined so that short curves correspond to admissible motions for the system. With convergence guarantee~\cite{AGHF2019}, we find this admissible motion by solving the AGHF. %Activation functions are defined to model the hybrid nature of legged locomotion. 
The method is able to automatically find the trajectory of CoM,  feet trajectory with contact positions and contact forces on even terrain all at once. Future work will include extension of the method to 3D models and automatically planning of schedule/timing of contact sequence.

\bibliographystyle{plain}
\bibliography{IEEEabrv,reference}

\end{document}